# ADMISSIBLE PREDICTIVE DENSITY ESTIMATION[1]

By Lawrence D. Brown, Edward I. George and Xinyi Xu

*University of Pennsylvania, University of Pennsylvania
and The Ohio State University*

Let $X|\mu \sim N_p(\mu, v_x I)$ and $Y|\mu \sim N_p(\mu, v_y I)$ be independent $p$-dimensional multivariate normal vectors with common unknown mean $\mu$. Based on observing $X = x$, we consider the problem of estimating the true predictive density $p(y|\mu)$ of $Y$ under expected Kullback–Leibler loss. Our focus here is the characterization of admissible procedures for this problem. We show that the class of all generalized Bayes rules is a complete class, and that the easily interpretable conditions of Brown and Hwang [*Statistical Decision Theory and Related Topics* (1982) *III* 205–230] are sufficient for a formal Bayes rule to be admissible.

**1. Introduction.** Let $X|\mu \sim N_p(\mu, v_x I)$ and $Y|\mu \sim N_p(\mu, v_y I)$ be independent $p$-dimensional multivariate normal vectors with common unknown mean $\mu \in R^p$. We assume that $v_x > 0$ and $v_y > 0$ are known. We let $p(x|\mu)$ and $p(y|\mu)$ denote the conditional densities of $X$ and $Y$, suppressing the dependence on $v_x$ and $v_y$ throughout.

Based on observing only $X = x$, we consider the problem of estimating the density $p(y|\mu)$ of $Y$. The natural action space $\mathcal{A}_0$ consists of all proper densities on $R^p$, that is

$$(1) \qquad \mathcal{A}_0 = \left\{ g : R^p \to R \text{ such that } g(y) \geq 0 \text{ and } \int g(y)\,dy = 1 \right\}.$$

For each observation $x \in R^p$, a (nonrandomized) decision procedure $\hat{p}(\cdot|x) : R^p \to \mathcal{A}_0$ chooses a $g \in \mathcal{A}_0$.

Received February 2006; revised May 2007.

[1]Supported by NSF Grants DMS-04-05716, DMS-06-05102 and DMS-07-07033.

We measure the goodness of fit of $g(y)$ to $p(y|\mu)$ by Kullback–Leibler (KL) loss

$$(2) \qquad L(\mu, g) = \begin{cases} \int p(y|\mu) \log \dfrac{p(y|\mu)}{g(y)} \, dy, & \text{if } g(y) > 0 \text{ a.e.,} \\ \infty, & \text{otherwise,} \end{cases}$$

and evaluate a procedure $\hat{p}(\cdot|x)$ by its risk function

$$(3) \qquad R_{\mathrm{KL}}(\mu, \hat{p}) = \int L(\mu, \hat{p}(\cdot|x)) p(x|\mu) \, dx.$$

For the comparison of two (nonrandomized) procedures, we say that $\hat{p}_1$ dominates $\hat{p}_2$ if $R_{\mathrm{KL}}(\mu, \hat{p}_1) \leq R_{\mathrm{KL}}(\mu, \hat{p}_2)$ for all $\mu$ and with strict inequality for some $\mu$. A procedure $\hat{p}(\cdot|x)$ is called admissible if it cannot be dominated.

Two widely-used methods to obtain predictive densities are "plug-in" rules and Bayes rules. A plug-in rule

$$(4) \qquad \hat{p}_{\hat{\mu}}(y|x) = p(y|\mu = \hat{\mu}(x))$$

simply substitutes an estimate $\hat{\mu}$ for $\mu$ in $p(y|\mu)$. In contrast, a Bayes rule integrates $\mu$ out with respect to a nonnegative and locally finite prior measure $M$ to obtain

$$(5) \qquad \hat{p}_M(y|x) = \frac{\int p(x|\mu) p(y|\mu) M(d\mu)}{\int p(x|\mu) M(d\mu)} = \int p(y|\mu) M(d\mu|x).$$

When writing an expression such as (5), we implicitly assume that the denominator in the middle expression is finite for all $x$, and hence all terms in (5) are finite for all $x$. We use the symbol $\pi$ to denote the density of $M$ when it exists, and will write either $\hat{p}_\pi$ or $\hat{p}_M$ in that case.

Aitchison (1975) showed that for proper $M$, $\hat{p}_M(y|x)$ minimizes the average KL risk

$$(6) \qquad B_{\mathrm{KL}}(M, \hat{p}) = \int R_{\mathrm{KL}}(\mu, \hat{p}) M(d\mu).$$

Aitchison also showed that the (formal) Bayes rule (5) under the uniform prior density $\pi_U(\mu) = 1$, namely $\hat{p}_{\pi_U}(y|x)$, dominates the plug-in rule $p(y|\hat{\mu}_{\mathrm{MLE}})$, which substitutes the maximum likelihood estimate $\hat{\mu}_{\mathrm{MLE}} = x$ for $\mu$. Indeed, as will be seen in Section 3, all the admissible procedures for multivariate normal density prediction under KL loss are Bayes rules in the sense of (5).

The constant risk Bayes rule $\hat{p}_{\pi_U}$ is best invariant, minimax, admissible when $p = 1$ [Murray (1977), Ng (1980) and Liang and Barron (2004)], and as we shall show in Section 3, admissible when $p = 2$. However, it is inadmissible when $p \geq 3$. This was first established by Komaki (2001) who showed that $\hat{p}_{\pi_U}$ is dominated by the Bayes rule under the (nonconstant) harmonic prior when $p \geq 3$. Liang (2002) further showed that $\hat{p}_{\pi_U}$ is dominated by proper Bayes rules under Strawderman priors when $p \geq 5$.

It is interesting to note the parallels between our predictive density estimation problem and the problem of estimating a multivariate normal mean under quadratic loss. Based on observing $Z|\mu \sim N_p(\mu, vI)$ with $v$ known, this latter problem is to estimate $\mu$ under quadratic risk

$$(7) \qquad R_Q^v(\mu, \hat{\mu}) = E_\mu \|\hat{\mu} - \mu\|^2,$$

where the dependence of $R_Q^v$ on $v$ is indicated by the superscript $v$. Here the maximum likelihood estimator $\hat{\mu}_{\text{MLE}}$, which is best invariant, minimax and admissible when $p = 1$ or 2, is dominated by the Bayes rule $\hat{\mu}_\pi = \int \mu \pi(\mu|x) \, d\mu$ under the harmonic prior when $p \geq 3$ [Stein (1974)] and by the proper Bayes rule under the Strawderman priors when $p \geq 5$ [Strawderman (1971)]. Note that in the KL risk problem $\hat{p}_{\pi_U}(y|x)$, rather than $\hat{p}(y|\hat{\mu}_{\text{MLE}})$, plays the same role as $\hat{\mu}_{\text{MLE}}$ in the quadratic risk problem. Recall that $\hat{\mu}_{\text{MLE}}$ can also be motivated as the Bayes rule under $\pi_U(\mu) = 1$ in the quadratic risk problem.

George, Liang and Xu (2006) recently drew out these parallels between the KL risk and quadratic risk problems, and found that they could be explained by connections between unbiased estimates of risk. These connections were shown to yield analogous sufficient conditions for the minimaxity of Bayes rules in both problems. In this paper, we establish further parallels concerning the characterization of admissibility in both problems. As proper Bayes rules are easily shown to be admissible in the KL setting, see Section 4.8.1 in Berger (1985), our focus will be on improper $\pi$ under which $\hat{p}_\pi(y|x)$ is sometimes more precisely called a formal or generalized Bayes rule. In Section 3, we establish sufficient conditions for the admissibility of Bayes rules $\hat{p}_\pi(y|x)$ under KL loss, conditions analogous to those of Brown (1971) and Brown and Hwang (1982). In Section 3, we prove that all admissible procedures for the KL risk problems are Bayes rules, a direct parallel of the complete class theorem of Brown (1971) for quadratic risk.

It might be of interest to note that when $v_y \to 0$, $p(y|\mu)$ degenerates to a point mass $I\{y = \mu\}$ and that by (5),

$$\hat{p}_\pi(y|x) = \int p(y|\mu)\pi(\mu|x) \, d\mu \to \pi(y|x).$$

Therefore, the limiting KL risk of a Bayes rule $\hat{p}_\pi$ is

$$\lim_{v_y \to 0} R_{\text{KL}}(\mu, \hat{p}_\pi) = E_\mu \left[ I\{y = \mu\} \log \frac{I\{y = \mu\}}{\pi(y|X)} \right] = -E_\mu \log \pi(\mu|X),$$

where the right-hand side can be viewed as the KL risk for "estimating a point mass at $\mu$" by a posterior density. Thus, our setup can provide a decision theoretic framework for evaluating a prior by the extent to which $E_\mu \log \pi(\mu|X)$ is large for all $\mu$.

**2. Sufficient conditions for admissibility.** For $Z \sim N_p(\mu, I)$, Brown (1971) and Brown and Hwang (1982) developed general sufficient conditions for the admissibility of formal Bayes rules for the quadratic risk problem. To utilize their results and obtain analogous sufficient conditions for the KL risk problem, we first establish a relationship between KL risk and quadratic risk. In this section, we assume that the prior measure $M$ has a density $\pi$ and that $R_{\mathrm{KL}}(\mu, \hat{p}_\pi) < \infty$ for all $\mu \in R^p$. Let

$$(8) \qquad m_\pi(z; v) = \int p(z|\mu)\pi(\mu)\,d\mu$$

be the marginal density of $Z \sim N_p(\mu, vI)$ under $\pi$.

THEOREM 1. *Let $\pi$ be a prior density on $\mu$ such that $m_\pi(z; v_x)$ is finite for all $z$. Then*

$$
\begin{aligned}
(9) \qquad & R_{\mathrm{KL}}(\mu, \hat{p}_{\pi_U}) - R_{\mathrm{KL}}(\mu, \hat{p}_\pi) \\
& = \frac{1}{2} \int_{v_w}^{v_x} \frac{1}{v^2} [R_Q^v(\mu, \hat{\mu}_{\mathrm{MLE}}) - R_Q^v(\mu, \hat{\mu}_\pi)]\,dv
\end{aligned}
$$

*where $v_w = v_x v_y/(v_x + v_y) < v_x$.*

PROOF. Let $m_\pi(w; v_w)$ denote the marginal density under $\pi$ of

$$(10) \qquad W = \frac{v_y X + v_x Y}{v_x + v_y} \sim N_p(\mu, v_w I).$$

By Lemmas 2 and 3 of George, Liang and Xu (2006),

$$
\begin{aligned}
(11) \qquad & R_{\mathrm{KL}}(\mu, \hat{p}_{\pi_U}) - R_{\mathrm{KL}}(\mu, \hat{p}_\pi) \\
& = E_{\mu, v_w} \log m_\pi(W; v_w) - E_{\mu, v_x} \log m_\pi(X; v_x),
\end{aligned}
$$

$m_\pi(z; v)$ is finite for any $v_w \le v \le v_x$, and

$$(12) \qquad \frac{\partial}{\partial v} E_{\mu,v} \log m_\pi(Z; v) = E_{\mu,v}\left(2\frac{\nabla^2 \sqrt{m_\pi(Z; v)}}{\sqrt{m_\pi(Z; v)}}\right),$$

where $\nabla^2 g(z) = \sum \frac{\partial^2}{\partial z_i^2} g(z)$, and $E_{\mu,v}(\cdot)$ stands for expectation with respect to the $N(\mu, vI)$ distribution. Furthermore, Stein (1974, 1981) showed that for the quadratic risk problem

$$(13) \qquad R_Q^v(\mu, \hat{\mu}_{\mathrm{MLE}}) - R_Q^v(\mu, \hat{\mu}_\pi) = -4v^2 E_{\mu,v}\frac{\nabla^2 \sqrt{m_\pi(Z; v)}}{\sqrt{m_\pi(Z; v)}}.$$

Combining (11), (12) and (13), the lemma follows.  □

Now let $B_{\mathrm{KL}}(\pi, \hat{p}) = \int R_{\mathrm{KL}}(\mu, \hat{p})\pi(\mu)\,d\mu$ and $B_Q^v(\pi, \hat{\mu}) = \int R_Q^v(\mu, \hat{\mu}) \times \pi(\mu)\,d\mu$ be the average KL and quadratic risks over $\pi$. The following relationship between the average KL risk difference and the average quadratic risk difference of Bayes rules follows from (9) and averaging over a prior $\pi_n$ that satisfies $\int_{R^p} \pi_n(\mu)\,d\mu < \infty$.

COROLLARY 1. *Let $\pi$ and $\pi_n$ be priors on $\mu$ such that $m_\pi(z; v_x)$ and $m_{\pi_n}(z; v_x)$ are finite for all $z$. Furthermore, assume $\pi_n$ satisfies $\int_{R^p} \pi_n(\mu)\,d\mu < \infty$. Then*

$$
\begin{aligned}
(14) \quad & B_{\mathrm{KL}}(\pi_n, \hat{p}_\pi) - B_{\mathrm{KL}}(\pi_n, \hat{p}_{\pi_n}) \\
& = \frac{1}{2} \int_{v_w}^{v_x} \frac{1}{v^2} [B_Q^v(\pi_n, \hat{\mu}_\pi) - B_Q^v(\pi_n, \hat{\mu}_{\pi_n})]\,dv.
\end{aligned}
$$

Corollary 1 enables us to extend the approach of Brown and Hwang (1982) to establish conditions for the admissibility of formal Bayes rules in the KL risk problem. As in Brown and Hwang (1982), we use Blyth's method which can be extended to any statistical estimation problem with a strictly convex loss function [Brown (1971)].

LEMMA 1. *Let $\hat{p}$ be such that $R_{\mathrm{KL}}(\mu, \hat{p}) < \infty$ for all $\mu \in R^p$. If there exists a sequence of densities $\{\pi_n\}$ such that $\int_{R^p} \pi_n(\mu)\,d\mu < \infty$, $\int_{\|\mu\| \leq 1} \pi_n(\mu)\,d\mu \geq c$ for some positive constant $c$, and*

$$
(15) \qquad B_{\mathrm{KL}}(\pi_n, \hat{p}) - B_{\mathrm{KL}}(\pi_n, \hat{p}_{\pi_n}) \to 0
$$

*then $\hat{p}$ is admissible.*

PROOF. Suppose $\hat{p}$ is not admissible. Then there is a $\hat{p}'$ such that $R_{\mathrm{KL}}(\mu, \hat{p}') \leq R_{\mathrm{KL}}(\mu, \hat{p})$ with strict inequality for some $\mu$. Let $\hat{p}'' = (\hat{p} + \hat{p}')/2$. Thus

$$
\begin{aligned}
& R_{\mathrm{KL}}(\mu, \hat{p}'') \\
& = \int\int p(x|\mu)p(y|\mu)\left[\log \frac{p(y|\mu)}{\hat{p}''(y|x)}\right]dx\,dy \\
& = \int\int p(x|\mu)p(y|\mu)\left[\log p(y|\mu) - \log\left(\frac{1}{2}\hat{p}(y|x) + \frac{1}{2}\hat{p}'(y|x)\right)\right]dx\,dy \\
& < \int\int p(x|\mu)p(y|\mu)\left[\log p(y|\mu) - \frac{1}{2}(\log\hat{p}(y|x) + \log\hat{p}'(y|x))\right]dx\,dy \\
& = \frac{1}{2}(R_{\mathrm{KL}}(\mu, \hat{p}) + R_{\mathrm{KL}}(\mu, \hat{p}')) \leq R_{\mathrm{KL}}(\mu, \hat{p}).
\end{aligned}
$$

Since $R_{\mathrm{KL}}(\mu, \hat{p})$ and $R_{\mathrm{KL}}(\mu, \hat{p}'')$ are both continuous in $\mu$, there exists an $\varepsilon > 0$ such that for all $\mu \in \{\mu : \|\mu\| \leq 1\}$,

$$
R_{\mathrm{KL}}(\mu, \hat{p}) - R_{\mathrm{KL}}(\mu, \hat{p}'') \geq \varepsilon > 0.
$$

Therefore, we have

$$B_{\mathrm{KL}}(\pi_n, \hat{p}) - B_{\mathrm{KL}}(\pi_n, \hat{p}_{\pi_n}) \geq B_{\mathrm{KL}}(\pi_n, \hat{p}) - B_{\mathrm{KL}}(\pi_n, \hat{p}'') \geq \varepsilon \cdot c > 0$$

which contradicts (15). The admissibility of $\hat{p}$ follows. $\square$

We assume without loss of generality that the coordinate system is chosen so that $\int_{\|\mu\| \leq 1} \pi(\mu) \, d\mu \geq c$ for some positive constant $c$. Using Lemma 1, we extend the approach of Brown and Hwang (1982) to obtain the following.

THEOREM 2. *A formal Bayes rule $\hat{p}_\pi$ is admissible under KL loss if for every $v \in [v_w, v_x]$, the improper $\pi$ satisfies both:*

(i) *the growth condition:*

$$\tag{16} \int_{R^p - S} \frac{\pi(\mu)}{\|\mu\|^2 \log^2(\|\mu\| \vee 2)} \, d\mu < \infty,$$

*where $S = \{\mu : \|\mu\| \leq 1\}$ and $a \vee b = \max\{a, b\}$, and*

(ii) *the asymptotic flatness condition:*

$$\tag{17} \int \int \pi(\mu) \left\| \frac{m_{\nabla \pi}(z; v)}{m_\pi(z; v)} - \frac{\nabla \pi}{\pi} \right\|^2 p(z | \mu) \, d\mu \, dz < \infty.$$

PROOF. For $v = 1$, Brown and Hwang showed that when the prior density $\pi$ satisfies the growth condition (16) and the asymptotic flatness condition (17), there exists a sequence of densities $\{\pi_n\}$ such that $\int_{\|\mu\| \leq 1} \pi_n(\mu) \, d\mu = \int_{\|\mu\| \leq 1} \pi(\mu) \, d\mu \geq c$ and that $B_Q^{v=1}(\pi_n, \hat{\mu}) - B_Q^{v=1}(\pi_n, \hat{\mu}_{\pi_n}) \to 0$. Furthermore, they showed that an explicit construction of such a sequence $\{\pi_n\}$ is obtained by defining

$$\tag{18} j_n(\mu) = \begin{cases} 1, & \|\mu\| \leq 1, \\ 1 - \dfrac{\log(\|\mu\|)}{\log(n)}, & 1 \leq \|\mu\| \leq n, \\ 0, & \|\mu\| \geq n, \end{cases}$$

for $n = 2, 3, \ldots,$ and letting

$$\tag{19} \pi_n(\mu) = j_n^2(\mu) \pi(\mu).$$

It is straightforward to show that the above construction also works for general $v$. That is, for any $v$, if $\pi$ satisfies conditions (16) and (17), then for the sequence $\{\pi_n\}$ obtained by (18) and (19), $\Delta_{n,v} \equiv B_Q^v(\pi_n, \hat{\mu}_\pi) - B_Q^v(\pi_n, \hat{\mu}_{\pi_n}) \to 0$. It thus follows that if $\pi$ satisfies conditions (16) and (17) for every $v \in [v_w, v_x]$, then by Corollary 1 and by the continuity in $v$ of $\Delta_{n,v}$,

$$\tag{20} B_{\mathrm{KL}}(\pi_n, \hat{p}_\pi) - B_{\mathrm{KL}}(\pi_n, \hat{p}_{\pi_n}) = \frac{1}{2} \int_{v_w}^{v_x} \frac{1}{v^2} \Delta_{n,v} \, dv \to 0.$$

That $\hat{p}_\pi$ is admissible now follows immediately from Lemma 1. $\square$

EXAMPLE 1 (Uniform prior). Let $\pi(\mu) = 1$ for any $\mu$, then $\nabla \pi = 0$. In this case, the conditions of Theorem 2 are easy to verify when $p = 1$ or 2. Therefore, the formal Bayes rule $\hat{p}_{\pi_U}$ is admissible when $p = 1$ or 2.

It was pointed out in Brown and Hwang (1982) that if

$$\pi(\mu) \leq \|\mu\|^{2-p},$$

(21)
$$\frac{\nabla \pi(\mu)}{\pi(\mu)} = o(\|\mu\|^{-1}) \quad \text{and} \quad \left| \frac{\partial^2 \pi(\mu)}{\partial \mu_i \partial \mu_j} \right| = o(\|\mu\|^{-2}),$$

hen (16) is easy to check and (17) can be verified with some difficulty [extending Lemma 3.4.1 of Brown (1971)]. Hence, by Theorem 2, the corresponding $\hat{p}_\pi$ is admissible under KL loss.

EXAMPLE 2 (*Harmonic prior*). Let $\pi_H(\mu) = \|\mu\|^{-(p-2)}$ for $p \geq 3$. Because this prior satisfies (21), the formal Bayes rule $\hat{p}_{\pi_H}$ is admissible when $p \geq 3$.

The following corollary is similarly a straightforward extension from Brown and Hwang (1982). It replaces condition (17) of Theorem 2 with a condition that is slightly less general, but more transparent and easier to verify.

COROLLARY 2. *If an improper density $\pi$ satisfies* (16) *and*

$$\int \frac{\|\nabla \pi(\mu)\|^2}{\pi(\mu)} \, d\mu < \infty,$$

(22)

*then the formal Bayes rule $\hat{p}_\pi$ is admissible under KL loss.*

Finally, it was also pointed out in Brown and Hwang (1982) that if

(23)    $\pi(\mu) \leq \|\mu\|^{2-p-\varepsilon}$    for some $\varepsilon > 0$    and    $\dfrac{\nabla \pi(\mu)}{\pi(\mu)} = o(\|\mu\|^{-1})$,

then (16) and (22) are easy to check. Hence, by Corollary 2, the corresponding $\hat{p}_\pi$ is admissible under KL loss.

There have been a few treatments of related problems yielding admissibility results in the same spirit as the above. In particular, Eaton (1982) formulates a prediction problem similar to the above, but under an integrated quadratic ($\mathcal{L}_2$) loss function, rather than our KL loss. Gatsonis (1984) discusses a related problem of estimating an unknown prior under this quadratic loss. Gatsonis proves an admissibility result in his setting for the Bayes procedure for the uniform prior. Gatsonis' methods do not easily apply to problems involving Bayes procedures for (generalized) priors other than the uniform prior. Eaton [(1992), Section 6] considers a prediction problem like ours, but with a different type of loss function. This loss function is bounded, and leads to a problem that is "quadratically regular" in a sense of that paper. For such quadratically regular problems the results of Eaton (1992) show admissibility for a specified class of prior measures. It is shown in Theorem 5.2 of Eaton et al. (2007) that Eaton's class of prior measures contains most of the densities covered by our Theorem 2, and vice-versa.

## 3. A complete class theorem.

We now turn to establishing that all (generalized) Bayes rules form a complete class for the KL loss problem. In Section 3.1, we begin by first establishing properties of some modified action spaces and the KL loss function. We then make use of these properties in Section 3.2 where we prove our main complete class results.

3.1. *Preliminary lemmas.* Because the true density $p(y|\mu)$ is bounded by a constant $C = (2\pi v_y)^{-p/2}$ for any $\mu$, it will eventually be useful to restrict attention to bounded density estimates. Let

$$(24) \qquad \mathcal{A} = \left\{ g : R^p \to R \text{ such that } 0 \leq g(y) \leq C \text{ a.e. and } \int g(y)\, dy = 1 \right\}.$$

Obviously, $\mathcal{A}$ is a subset of the action space $\mathcal{A}_0$ that is defined in (1).

The following lemma, which is proved in the Appendix, shows that no admissible actions are lost by restricting the action space to $\mathcal{A}$.

LEMMA 2. *Suppose $g_0(\cdot) \in \mathcal{A}_0$. If $g_0 \notin \mathcal{A}$, that is, $g_0 > C$ on a set $S \subset R^p$ with positive measure, then there exists a $g \in \mathcal{A}$ that dominates $g_0$ in the sense that $L(\mu, g_0) > L(\mu, g)$ for all $\mu$.*

It will also be useful to consider extending $\mathcal{A}$ to its closure

$$(25) \qquad \mathcal{A}^* = \left\{ g : R^p \to R \text{ such that } 0 \leq g(y) \leq C \text{ a.e. and } \int g(y)\, dy \leq 1 \right\},$$

and then to make use of the topological properties of $\mathcal{A}^*$. Because $\mathcal{A}^*$ is a subset of the Banach space $\mathcal{L}_\infty$, we will consider the topology on $\mathcal{A}^*$ induced by the weak* topology on $\mathcal{L}_\infty$. Under this weak* topology, a sequence $\{g_i\} \in \mathcal{A}^*$ converges to a $g \in \mathcal{A}^*$ if

$$(26) \qquad \int f(y) g_i(y)\, dy \to \int f(y) g(y)\, dy \qquad \forall f \in \mathcal{L}_1.$$

We will eventually make use of the following properties of $\mathcal{A}^*$ under the weak* topology.

LEMMA 3. *Define the action space $\mathcal{A}^*$ as in (25), then:*

(i) $\mathcal{A}^*$ is weak* compact.

(ii) The weak* topology on $\mathcal{A}^*$ is metrizable by

$$(27) \qquad \rho(g, h) = \sum_{k=1}^{\infty} 2^{-k} \left| \int [g(y) - h(y)] f_k(y) \, dy \right| \qquad \text{for any } g, h \in \mathcal{A}^*,$$

where $\{f_k, k = 1, 2, \ldots\}$ is a countable dense subset of $\mathcal{L}_1$. And $\mathcal{A}^*$ is separable and second countable under this metric $(27)$.

(iii) Suppose $g^*(\cdot) \in \mathcal{A}^*$. If $g^* \notin \mathcal{A}$, then there exists a $g \in \mathcal{A}$ that dominates $g^*$ in the sense that $L(\mu, g^*) > L(\mu, g)$ for all $\mu$. Thus, the extension from $\mathcal{A}$ to $\mathcal{A}^*$ does not incur any new admissible actions.

Finally, we also need to make use of the following properties of the Kullback–Leibler loss function.

LEMMA 4. For the KL loss function $L(\mu, \cdot)$ in $(2)$:

(i) $L(\mu, \cdot)$ is lower semi-continuous on $\mathcal{A}^*$, that is, if $\{g_i\}, g \in \mathcal{A}^*$ and $g_i \to g \in \mathcal{A}^*$ weak*, then

$$(28) \qquad \liminf_{i \to \infty} L(\mu, g_i) \geq L(\mu, g) \qquad \forall \mu \in R^p;$$

(ii) $L(\mu, \cdot)$ is strictly convex on

$$(29) \qquad \mathcal{A}_+^* = \{g : g \in \mathcal{A}^* \text{ and } L(\mu, g) < \infty \text{ for } \forall \mu\}$$

for any $\mu \in R^p$.

3.2. The main theorems. Having established Lemmas 2, 3 and 4 in Section 3.1, we are now ready to prove that all admissible procedures for the normal density prediction problem under KL loss are (generalized) Bayes rules. This proof consists of three steps:

(i) All the admissible procedures are nonrandomized.

(ii) For any admissible procedure $\hat{p}(\cdot|x)$, there exists a sequence of priors $M_i(\mu)$ such that $\hat{p}_{M_i}(\cdot|x) \to \hat{p}(\cdot|x)$ for almost every $x$ under the weak* topology $(26)$.

(iii) We can find a subsequence $\{M_{i'}\}$ and a limit prior $M$ such that $\hat{p}_{M_{i'}}(\cdot|x) \to \hat{p}_M(\cdot|x)$ weak* for almost every $x$. Therefore, $\hat{p}(\cdot|x) = \hat{p}_M(\cdot|x)$ for a.e. $x$, that is, $\hat{p}(\cdot|x)$ is a (generalized) Bayes rule.

THEOREM 3. All nonrandomized procedures form a complete class.

PROOF. Let $\delta : R^p \to P(\mathcal{A}_0)$ be an admissible and randomized procedure, where $P(\mathcal{A}_0)$ denotes the space of probability distributions over $\mathcal{A}_0$. We first prove that $\delta(x) \in P(\mathcal{A}) \subset P(\mathcal{A}^*)$ for a.e. $x$. Suppose there exists a set $K$ such that $K$ has positive measure and for each $x \in K$, $\delta(\cdot|x) = \hat{p}(\cdot|x) \notin P(\mathcal{A}^*)$ with a positive probability. Then by Lemma 2, we can find $g_x \in P(\mathcal{A}^*)$ that satisfies $L(\mu, g_x) < L(\mu, \hat{p}(\cdot|x))$ for all $\mu$, and therefore $\delta$ is dominated by the decision rule $\tilde{\delta}$ that substitutes $g_x$ for $\hat{p}(\cdot|x)$. This contradicts the admissibility of $\delta$.

Now let $\hat{p}^*(y|x) = E^{\delta(\cdot|x)}(g(y))$. It can be seen that $\hat{p}^*(y|x) \in \mathcal{A}$ since $\delta(x) \in P(\mathcal{A}) \subset P(\mathcal{A}^*)$ for a.e. $x$. By Lemma 4(ii) and Jensen's inequality,

$$(30) \qquad L(\mu, \hat{p}^*(y|x)) \le E^{\delta(\cdot|x)}(L(\mu, \hat{p}(y))) = L(\mu, \delta(y|x)) \qquad \forall \mu.$$

Furthermore, strict inequality holds in (30) unless either $\delta(\cdot|x)$ is nonrandomized with probability 1 or $L(\mu, \delta(y|x)) = \infty$, which implies that $\delta$ can be dominated by a finite-risk nonrandomized procedure. Therefore, it contradicts that $\delta$ is admissible and randomized. It then follows that the nonrandomized procedures are a complete class. □

Theorem 3 shows that we can restrict attention to nonrandomized procedures $\hat{p}(\cdot|x)$. Next we prove that for a.e. $x$, all admissible procedures are limits of Bayes rules (5). Since the Bayes rules are also nonrandomized, this convergence can be evaluated with respect to the weak* topology for each $x$.

THEOREM 4. *For any admissible procedure $\hat{p}(\cdot|x)$, there exists a sequence of priors $\{M_i\}$ supported on finite sets such that $\hat{p}_{M_i}(\cdot|x) \to \hat{p}(\cdot|x)$ weak* for a.e. $x$ under the topology* (26).

PROOF. This is essentially Theorem 4A.12 of Brown (1986). There are some minor differences between the formulations there and here which we now note in order to clarify how that Theorem 4A.12 yields the current Theorem 4. The principal difference is that the action space $\mathcal{A}^*$ in Brown (1986) was assumed to be Euclidean whereas here it is merely compact, separable, and metrizable. Because the space $\mathcal{A}^*$, here is compact, the one-point compactification $\{i\}$ introduced in Brown (1986) is not needed. This simplifies the proof of Proposition 4A.11 there, which in our context becomes Theorem 3. The remainder of the proof proceeds as discussed in the text of the proof of Theorem 4A.12. □

Theorem 4 establishes that any admissible procedure $\hat{p}(\cdot|x)$ is a limit of Bayes rules for a.e. $x$. To prove $\hat{p}(\cdot|x)$ itself is also a (generalized) Bayes rule, we need to find a (possibly improper) prior $M$ such that $\hat{p}_M(\cdot|x) = \hat{p}(\cdot|x)$ for a.e. $x$.

THEOREM 5. *The set of all generalized Bayes procedures is a complete class of procedures.*

PROOF. Suppose $\hat{p}(\cdot|x)$ is an admissible procedure. Then by Theorem 4, there exists a sequence of measures $M_i$ supported on finite sets such that $\hat{p}_{M_i}(\cdot|x) \to \hat{p}(\cdot|x)$ for a.e. $x$ under the weak* topology (26).

Let
$$r_i = \int_{\|x\| \leq 1} \int p(x|\mu) M_i(d\mu)\, dx,$$

then $r_i > 0$ since $p(x|\mu) > 0$ for all $x$ and $\mu$. Thus we can define a new sequence of measures $M_i'$ by $M_i' = M_i/r_i$. It is easy to check that $\hat{p}_{M_i'} = \hat{p}_{M_i} \to \hat{p}$ weak* a.e. and that

$$(31) \qquad \int_{\|x\| \leq 1} \int p(x|\mu) M_i'(d\mu)\, dx = 1.$$

By 2.16(iv) of Brown (1986), there exists a finite limiting measure $M$ such that $M_i' \to M$.

Let $S$ be the biggest convex set that satisfies
$$\liminf_{i \to \infty} \sup_{x \in S} \int p(x|\mu) M_i'(d\mu) < \infty.$$

[The existence of $S$ follows from (31).] Then by Theorem 2.17 in Brown (1986), for any $x$ in the interior of $S$,

$$(32) \qquad \int p(x|\mu) M_i'(d\mu) \to \int p(x|\mu) M(d\mu) \qquad \text{as } i \to \infty.$$

In fact, we can prove that the closure $\bar{S} = R^p$. [Otherwise its complement $\bar{S}^c$ has positive measure and at every $x \in \bar{S}^c$, $\liminf_{i \to \infty} \int p(x|\mu) M_i'(d\mu) = \infty$. Therefore,

$$\lim_{i \to \infty} \int_{\|y\| \leq 1} \hat{p}_{M_i'}(y|x)\, dy = \lim_{i \to \infty} \int_{\|y\| \leq 1} \frac{\int p(x|\mu) p(y|\mu) M_i'(d\mu)}{\int p(x|\mu) M_i'(d\mu)}\, dy$$
$$\leq (2\pi v_x)^{-p/2} \lim_{i \to \infty} \frac{\int M_i'(d\mu) \int_{\|y\| \leq 1} p(y|\mu)\, dy}{\int p(x|\mu) M_i'(d\mu)}$$
$$= 0,$$

which implies $\int_{\|y\| \leq 1} \hat{p}(y|x)\, dy = 0$ and thus $R_{\mathrm{KL}}(\mu, \hat{p}) = \infty$. This would contradict the assumed admissibility of $\hat{p}$.] Hence, (32) holds for a.e. $x$.

Furthermore, by the dominated convergence, for a.e. $x$ and $y$,

$$(33) \qquad \int p(x|\mu) p(y|\mu) M_i'(d\mu) \to \int p(x|\mu) p(y|\mu) M(d\mu).$$

Combining (32) and (33), we obtain

$$\hat{p}_{M_i'} = \frac{\int p(x|\mu) p(y|\mu) M_i'(d\mu)}{\int p(x|\mu) M_i'(d\mu)} \to \frac{\int p(x|\mu) p(y|\mu) M(d\mu)}{\int p(x|\mu) M(d\mu)} = \hat{p}_M(y|x)$$

for a.e. $x$ and $y$, so $\hat{p}_{M_i'}$ also converges to $\hat{p}_M(y|x)$ under the weak* topology. Therefore, $\hat{p} = \hat{p}_M$ is a generalized Bayes procedure.  $\square$

## APPENDIX

In this appendix, we provide the proofs of Lemmas 2, 3 and 4 from Section 3.1.

PROOF OF LEMMA 2.   (i) Suppose $g_0 = 0$ on a set with positive measure. Then by definition $L(\mu, g_0) = \infty$ for any $\mu$. So any $g \in \mathcal{A}$ with finite risk dominates it and thus $g_0$ is inadmissible.

(ii) Suppose $g_0 > 0$ almost everywhere. If $g_0 \geq C$ on a set $S$ with Lebesgue measure $\nu(S) > 0$, then a $g$ can be constructed by truncating $g_0$ on $S$ and lifting it in the other areas. Notice that $\int_{S^c} g_0(y)\, dy > 0$, so we can define

$$(34) \qquad c = \frac{1 - C\nu(S)}{\int_{S^c} g_0(y)\, dy},$$

where $S^c$ is the complement of $S$. It is easy to check $c > 1$. Let

$$(35) \qquad g(y) = \begin{cases} cg_0, & y \in S^c, \\ C, & y \in S. \end{cases}$$

Obviously, $g \in \mathcal{A}$. For any $\mu$, the difference between the loss functions of $g_0$ and $g$ is

$$L(\mu, g_0) - L(\mu, g)$$

$$= \int p(y|\mu) \log g(y)\, dy - \int p(y|\mu) \log g_0(y)\, dy$$

$$= \int_S p(y|\mu) \log C\, dy + \int_{S^c} p(y|\mu) \log(cg_0(y))\, dy - \int p(y|\mu) \log g_0(y)\, dy$$

$$= \int_S p(y|\mu) \log C\, dy + \log c \int_{S^c} p(y|\mu)\, dy$$

$$\quad + \int_{S^c} p(y|\mu) \log g_0(y)\, dy - \int p(y|\mu) \log g_0(y)\, dy$$

$$= \int_S p(y|\mu) \log C\, dy + \log c \int_{S^c} p(y|\mu)\, dy - \int_S p(y|\mu) \log g_0(y)\, dy$$

$$= \log c - \int_S p(y|\mu) \log \frac{cg_0(y)}{C}\, dy$$

$$\geq \log c - \log \int_S p(y|\mu) \frac{cg_0(y)}{C}\, dy \qquad \text{(Jensen's inequality)}$$

$$\geq \log c - \log \int_S cg_0(y)\, dy$$

$$> 0.$$

The last strict inequality holds because $\int_S g_0(y)\, dy = 1 - \int_{S^c} g_0(y)\, dy < 1$. Therefore, $g$ dominates $g_0$.   $\square$

PROOF OF LEMMA 3. (i) By the Banach–Alaoglu theorem, the $\mathcal{L}_1$ unit ball $\{g : R^p \to R | \int g(y)\,dy \leq 1\}$ is weak* compact. Also, it is easy to check that the bounded set $\{g : R^p \to R | 0 \leq g(y) \leq C\}$ is closed and thus compact. So their intersection $\mathcal{A}^*$ is compact.

(ii) Because $\mathcal{L}_1$ a separable normed space, the weak* topology on the closed ball of its dual space $\mathcal{L}_\infty$ can be metrized by (27). And since every compact metric space is separable and second countable, (ii) follows immediately from (i).

(iii) Suppose $g^* \in \mathcal{A}^*$ but $g^* \notin \mathcal{A}$, then $\int g^*(y)\,dy < 1$. If $\int g^*(y)\,dy = 0$, its loss function $L(\mu, g^*) = \infty$ for any $\mu$ and thus $g^*$ is inadmissible. Otherwise let $g' = g^* / \int g^*(y)\,dy$, then $\int g'(y)\,dy = 1$ and it is easy to check that $g'$ dominates $g^*$. Truncate $g'$ as in (35) if necessary, and then it yields a $g \in \mathcal{A}$ that dominates $g'$ and therefore dominates $g^*$. $\square$

PROOF OF LEMMA 4. (i) Suppose $\{g_i\}$ is a sequence of functions in $\mathcal{A}^*$ and $g_i \to g \in \mathcal{A}^*$ under the weak* topology.

(a) We first consider the case where $g$ is bounded away from 0. To prove that $\liminf_{i\to\infty} L(\mu, g_i) \geq L(\mu, g)$ for all $\mu \in R^p$, we only need to show

$$
\begin{aligned}
& L(\mu, g) - \liminf_{i\to\infty} L(\mu, g_i) \\
(36) \quad &= \limsup_{i\to\infty} \int p(y|\mu) \log g_i(y)\,dy - \int p(y|\mu) \log g(y)\,dy \\
&\leq 0.
\end{aligned}
$$

If there exists a positive constant $\varepsilon_0$ such that $g > \varepsilon_0$ a.e., then $\frac{p(y|\mu)}{g(y)} \leq \frac{p(y|\mu)}{\varepsilon_0}$ is an $\mathcal{L}_1$ function. Therefore,

$$
\begin{aligned}
& \limsup_{i\to\infty} \int p(y|\mu) \log g_i(y)\,dy - \int p(y|\mu) \log g(y)\,dy \\
&= \limsup_{i\to\infty} \int p(y|\mu) \log \frac{g_i(y)}{g(y)}\,dy \\
(37) \quad &\leq \limsup_{i\to\infty} \int p(y|\mu) \left( \frac{g_i(y)}{g(y)} - 1 \right)\,dy \\
&= \limsup_{i\to\infty} \int \frac{p(y|\mu)}{g(y)} g_i(y)\,dy - 1 \\
&= 0,
\end{aligned}
$$

where the inequality follows from the fact that $\log x \leq x - 1$ for any $x > 0$. This proves that the lemma holds whenever $g$ is bounded away from 0.

(b) Let $N = \{y : g(y) = 0\}$. If $N$ has positive measure, then the assumption that $g_i \to g$ under the weak* topology implies that $g_i(y) \to 0$ in measure on $N$. Hence

$$\lim_{i \to \infty} \int p(y|\mu) \log g_i(y)\, dy = -\infty \tag{38}$$

by the bounded convergence theorem for convergence in measure.

(c) Now the final possibility is that $N$ has measure 0, but $g$ is not bounded away from 0. Then for any fixed $\varepsilon > 0$, let $L(\varepsilon) = \{y | g(y) \geq \varepsilon\}$. Thus,

$$
\begin{aligned}
\limsup_{i \to \infty} & \int p(y|\mu) \log g_i(y)\, dy \\
&= \limsup_{i \to \infty} \left[ \int_{L(\varepsilon)} p(y|\mu) \log g_i(y)\, dy + \int_{L^c(\varepsilon)} p(y|\mu) \log g_i(y)\, dy \right] \\
&\leq \int_{L(\varepsilon)} p(y|\mu) \log g(y)\, dy + \log C \int_{L^c(\varepsilon)} p(y|\mu)\, dy.
\end{aligned}
\tag{39}
$$

The above inequality follows from the truth of the lemma when $g$ is bounded away from 0 and the definition that $g_i \in \mathcal{A}^*$ satisfies $g_i \leq C$. Now let $\varepsilon \downarrow 0$, then $L(\varepsilon) \to R^p$ since $g > 0$ a.e. Therefore, by the bounded convergence theorem,

$$
\begin{aligned}
\limsup_{i \to \infty} \int p(y|\mu) \log g_i(y)\, dy &\leq \int p(y|\mu) \log g(y)\, dy + 0 \\
&= \int p(y|\mu) \log g(y)\, dy.
\end{aligned}
\tag{40}
$$

This proves (i) since

$$L(\mu, g_i) = \int p(y|\mu) \log p(y|\mu)\, dy - \int p(y|\mu) \log g_i(y)\, dy. \tag{41}$$

(ii) Suppose $g_1, g_2 \in \mathcal{A}_+^*$ and $g_\lambda(y|x) = \lambda g_1 + (1-\lambda) g_2$ with $0 < \lambda < 1$, then

$$
\begin{aligned}
L(\mu, g_\lambda) &= \int p(y|\mu) \log \frac{p(y|\mu)}{g_\lambda(y)}\, dy \\
&= \int p(y|\mu) \log p(y|\mu)\, dy - \int p(y|\mu) \log[\lambda g_1(y) + (1-\lambda) g_2(y)]\, dy \\
&< \int p(y|\mu) \log p(y|\mu)\, dy - \int p(y|\mu)[\lambda \log g_1(y) \\
&\qquad\qquad\qquad\qquad\qquad\qquad + (1-\lambda) \log g_2(y)]\, dy \\
&= \lambda \int p(y|\mu) \log \frac{p(y|\mu)}{g_1(y)}\, dy + (1-\lambda) \int p(y|\mu) \log \frac{p(y|\mu)}{g_2(y)}\, dy \\
&= \lambda L(\mu, g_1) + (1-\lambda) L(\mu, g_2),
\end{aligned}
$$

where the inequality follows from Jensen's inequality. Thus, the strict convexity of $L(\mu, \cdot)$ on $\mathcal{A}_+^*$ is verified. $\square$

# REFERENCES

AITCHISON, J. (1975). Goodness of prediction fit. *Biometrika* **62** 547–554. MR0391353

BERGER, J. O. (1985). *Statistical Decision Theory and Bayesian Analysis*, 2nd ed. Springer, New York. MR0804611

BROWN, L. D. (1971). Admissible estimators, recurrent diffusions, and insoluble boundary value problems. *Ann. Math. Statist.* **42** 855–903. MR0286209

BROWN, L. D. (1986). *Fundamentals of Statistical Exponential Families with Applications in Statistical Decision Theory.* IMS, Hayward, CA. MR0882001

BROWN, L. D. and HWANG, J. (1982). A unified admissibility proof. In *Statistical Decision Theory and Related Topics III* (S. S. Gupta and J. O. Berger, eds.) **1** 205–230. Academic Press, New York. MR0705290

EATON, M. L. (1982). A method for evaluating improper prior distributions. In *Statistical Decision Theory and Related Topics III* (S. S. Gupta and J. O. Berger, eds.) **1** 329–352. Academic Press, New York. MR0705296

EATON, M. L. (1992). A statistical diptych: Admissible inferences–recurrence of symmetric Markov chains. *Ann. Statist.* **20** 1147–1179. MR1186245

EATON, M. L., HOBERT, J. P., JONES, G. L. and LAI, W.-L. (2007). Evaluation of formal posterior distributions via Markov chain arguments. Preprint. Available at http://www.stat.ufl.edu/~jhobert/.

GATSONIS, C. A. (1984). Deriving posterior distributions for a location parameter: A decision theoretic approach. *Ann. Statist.* **12** 958–970. MR0751285

GEORGE, E. I., LIANG, F. and XU, X. (2006). Improved minimax prediction under Kullback–Leibler loss. *Ann. Statist.* **34** 78–91. MR2275235

KOMAKI, F. (2001). A shrinkage predictive distribution for multivariate normal observations. *Biometrika* **88** 859–864. MR1859415

LIANG, F. (2002). Exact minimax procedures for predictive density estimation and data compression. Ph.D. dissertation, Dept. Statistics, Yale Univ.

LIANG, F. and BARRON, A. (2004). Exact minimax strategies for predictive density estimation, data compression and model selection. *IEEE Trans. Inform. Theory* **50** 2708–2726. MR2096988

MURRAY, G. D. (1977). A note on the estimation of probability density functions. *Biometrika* **64** 150–152. MR0448690

NG, V. M. (1980). On the estimation of parametric density functions. *Biometrika* **67** 505–506. MR0581751

STEIN, C. (1974). Estimation of the mean of a multivariate normal distribution. In *Proceedings of the Prague Symposium on Asymptotic Statistics* (J. Hajek, ed.) 345–381. Univ. Karlova, Prague. MR0381062

STEIN, C. (1981). Estimation of a multivariate normal mean. *Ann. Statist.* **9** 1135–1151. MR0630098

STRAWDERMAN, W. E. (1971). Proper Bayes minimax estimators of the multivariate normal mean. *Ann. Math. Statist.* **42** 385–388. MR0397939

L. D. Brown
E. I. George
Statistics Department
The Wharton School
University of Pennsylvania
Philadelphia, Pennsylvania 19104-6340
USA
E-mail: lbrown@wharton.upenn.edu
        edgeorge@wharton.upenn.edu

X. Xu
Department of Statistics
The Ohio State University
Columbus, Ohio 43210-1247
USA
E-mail: xinyi@stat.ohio-state.edu